\documentclass[11pt,fleqn]{article}
\usepackage{amsfonts}
\newcommand{\ol}{\setlength{\itemsep}{0pt.}\begin{enumerate}}
\newcommand{\eol}{\end{enumerate}\setlength{\itemsep}{-\parsep}}
\newcommand{\ignore}[1]{}
\title{Inverse Conjecture for the Gowers norm is false}
\author{
Shachar Lovett
\thanks{
Faculty of Mathematics and Computer Science,
The Weizmann Institute of Science,
Rehovot 76100, Israel.
Research supported by
ISF grant 1300/05.}
\and
Roy Meshulam
\thanks{Department of Mathematics,
The Technion,
Haifa 32000, Israel.}
\and
Alex Samorodnitsky
\thanks{School of Engineering and Computer Science,
The Hebrew University of Jerusalem,
Jerusalem 91904, Israel.
Research partially supported by ISF
grant 039-7682 and by BSF grant 2006377.}
}

\begin{document}
\date{}
\maketitle


\newtheorem{THEOREM}{Theorem}[section]
\newenvironment{theorem}{\begin{THEOREM} \hspace{-.85em} {\bf :}
}%
                        {\end{THEOREM}}
\newtheorem{LEMMA}[THEOREM]{Lemma}
\newenvironment{lemma}{\begin{LEMMA} \hspace{-.85em} {\bf :} }%
                      {\end{LEMMA}}
\newtheorem{COROLLARY}[THEOREM]{Corollary}
\newenvironment{corollary}{\begin{COROLLARY} \hspace{-.85em} {\bf
:} }%
                          {\end{COROLLARY}}
\newtheorem{PROPOSITION}[THEOREM]{Proposition}
\newenvironment{proposition}{\begin{PROPOSITION} \hspace{-.85em}
{\bf :} }%
                            {\end{PROPOSITION}}
\newtheorem{DEFINITION}[THEOREM]{Definition}
\newenvironment{definition}{\begin{DEFINITION} \hspace{-.85em} {\bf
:} \rm}%
                            {\end{DEFINITION}}
\newtheorem{EXAMPLE}[THEOREM]{Example}
\newenvironment{example}{\begin{EXAMPLE} \hspace{-.85em} {\bf :}
\rm}%
                            {\end{EXAMPLE}}
\newtheorem{CONJECTURE}[THEOREM]{Conjecture}
\newenvironment{conjecture}{\begin{CONJECTURE} \hspace{-.85em}
{\bf :} \rm}%
                            {\end{CONJECTURE}}
\newtheorem{MAINCONJECTURE}[THEOREM]{Main Conjecture}
\newenvironment{mainconjecture}{\begin{MAINCONJECTURE} \hspace{-.85em}
{\bf :} \rm}%
                            {\end{MAINCONJECTURE}}
\newtheorem{PROBLEM}[THEOREM]{Problem}
\newenvironment{problem}{\begin{PROBLEM} \hspace{-.85em} {\bf :}
\rm}%
                            {\end{PROBLEM}}
\newtheorem{QUESTION}[THEOREM]{Question}
\newenvironment{question}{\begin{QUESTION} \hspace{-.85em} {\bf :}
\rm}%
                            {\end{QUESTION}}
\newtheorem{REMARK}[THEOREM]{Remark}
\newenvironment{remark}{\begin{REMARK} \hspace{-.85em} {\bf :}
\rm}%
                            {\end{REMARK}}

\newcommand{\thm}{\begin{theorem}}
\newcommand{\lem}{\begin{lemma}}
\newcommand{\pro}{\begin{proposition}}
\newcommand{\dfn}{\begin{definition}}
\newcommand{\rem}{\begin{remark}}
\newcommand{\xam}{\begin{example}}
\newcommand{\cnj}{\begin{conjecture}}
\newcommand{\mcnj}{\begin{mainconjecture}}
\newcommand{\prb}{\begin{problem}}
\newcommand{\que}{\begin{question}}
\newcommand{\cor}{\begin{corollary}}
\newcommand{\prf}{\noindent{\bf Proof:} }
\newcommand{\ethm}{\end{theorem}}
\newcommand{\elem}{\end{lemma}}
\newcommand{\epro}{\end{proposition}}
\newcommand{\edfn}{\bbox\end{definition}}
\newcommand{\erem}{\bbox\end{remark}}
\newcommand{\exam}{\bbox\end{example}}
\newcommand{\ecnj}{\bbox\end{conjecture}}
\newcommand{\emcnj}{\bbox\end{mainconjecture}}
\newcommand{\eprb}{\bbox\end{problem}}
\newcommand{\eque}{\bbox\end{question}}
\newcommand{\ecor}{\end{corollary}}
\newcommand{\eprf}{\bbox}
\newcommand{\beqn}{\begin{equation}}
\newcommand{\eeqn}{\end{equation}}
\newcommand{\blist}{\begin{itemize}}
\newcommand{\elist}{\end{itemize}}
\newcommand{\wbox}{\mbox{$\sqcap$\llap{$\sqcup$}}}
\newcommand{\bbox}{\vrule height7pt width4pt depth1pt}
\newcommand{\qed}{\bbox}
\def\sup{^}

\def\B{\{0,1\}}
\def\H{\{-1,1\}}

\def\S{S(n,w)}

\def\n{\lfloor \frac n2 \rfloor}

\def\Tp{Tchebyshef polynomial}
\def\Tps{TchebysDeto be the maximafine $A(n,d)$ l size of a code with
  distance $d$hef polynomials}
\newcommand{\rarrow}{\rightarrow}

\newcommand{\larrow}{\leftarrow}

\overfullrule=0pt
\def\setof#1{\lbrace #1 \rbrace}

\def \E{\mathbb E}
\def \R{\mathbb R}
\def \C{\mathbb C}
\def \Z{\mathbb Z}
\def \FF{\mathbb F}
\def \FN{{\FF}^N}

\def\<{\left<}
\def\>{\right>}
\def \({\left(}
\def \){\right)}
\def \[{\left[}
\def \]{\right]}

\def \e{\epsilon}
\def \O{\Omega}
\def \N{{\cal N}}
\def \l{{\lambda}}
\def \S{{\cal S}}
\def \F{{\cal F}}
\def \H{{\cal H}}
\def \1{{\bf 1}}
\def \T{{\tau}}
\def \b{{\bf b}}
\def \a{{\bf a}}

\begin{abstract}
Let $p$ be a fixed prime number, and $N$ be a large integer.
The 'Inverse Conjecture for the Gowers norm' states that if the "$d$-th Gowers norm" of a function $f:~\FF^N_p \to \FF$ is non-negligible, that is larger than a constant independent of $N$, then $f$ can be non-trivially approximated by a degree $d-1$ polynomial. The conjecture is known to hold for $d=2,3$ and for any prime $p$. In this paper we show the conjecture to be false for $p=2$ and for $d = 4$, by presenting an explicit function whose $4$-th Gowers norm is non-negligible, but whose correlation any polynomial of degree $3$ is exponentially small.

Essentially the same result (with different correlation bounds) was independently obtained by Green and Tao \cite{gt07}. Their analysis uses a modification of a Ramsey-type argument of Alon and Beigel \cite{ab} to show inapproximability of certain functions by low-degree polynomials.

We observe that a combination of our results with the argument of Alon and Beigel implies the inverse conjecture to be false for any prime $p$, for $d = p^2$.
\end{abstract}

\section{Introduction}
We consider multivariate functions over finite fields. The main question of interest here would be whether these functions can be non-trivially approximated by a low-degree polynomial.

Fix a prime number $p$. Let $\FF = \FF_p$ be the finite field with $p$ elements.
Let $\xi = e^{\frac{2\pi i}{p}}$ be the primitive $p$-th root of unity. Denote by $e(x)$ the exponential function taking $x \in \FF$ to $\xi^x \in \C$. For two functions $f,g:~\FF^N \to \FF$, let $\<f,g\> := \E_x e(f(x) - g(x))$.
\dfn
\label{dfn:approx}
A function $f$ is non-trivially approximable by a degree-$d$ polynomial if
$$
|\<f,g\>| > \epsilon
$$
for some polynomial $g$ of degree at most $d$ in $\FF[x_1...x_N]$.
\edfn

More precisely, in this definition we are looking at a sequence $f_N$ of functions and of approximating low-degree polynomials $g_N$ in $N$ variables, and let $N$ grow to infinity. In this paper, the remaining parameters, that is the field size $p$, the degree $d$ and the offset $\e$ are fixed, independent of $N$.

A counting argument shows that a generic function can not be approximated by a polynomial of low degree. The problems of showing a specific given function to have no non-trivial approximation and of constructing an explicit non-approximable function have been extensively investigated, since solutions to these problems have many applications in complexity (cf. discussion and references in \cite{ab,vw,bv}).

This paper studies a technical tool that measures distance from low-degree polynomials. This is the Gowers norm, introduced in \cite{gowers}. For a function $f:~\FF^N \rarrow \FF$ and a vector $y \in \FF^n$, we take $f_y$ to be the directional derivative of $f$ in direction $y$ by setting
$$
f_y(x) = f(x+y) - f(x)
$$
For a $k$-tuple of vectors $y_1...y_k$ we take the iterated derivative in these directions to be
$$
f_{y_1...y_k} = \(f_{y_1...y_{k-1}}\)_{y_k}
$$
It is easy to see that this definition does not depend on the ordering of $y_1...y_k$.

The $k$-th Gowers "norm" $\|f\|_{U^k}$ of $f$ is
$$
\(\E_{x,y_1...y_k} \[e\(f_{y_1...y_k}(x)\)\]\)^{1/2^k}
$$
More accurately, as shown in \cite{gowers}, this is indeed a norm of the associated complex-valued function $e(f)$ (for $k \ge 2$).

It is easy to see that $\|f\|_{U^{d+1}}$ is $1$ iff $f$ is a polynomial of degree at most $d$. This is just another way of saying that all order-$(d+1)$ iterative derivatives of $f$ are zero if and only if $f$ is a polynomial of degree at most $d$. It is also possible to see \cite{gt} that $|\<f,g\>| > \e$ for $g$ of degree at most $d$, implies $\|f\|_{U^{d+1}} > \e$. That is to say, if $f$ is non-trivially close to a degree-$d$ polynomial, this can be detectable via an appropriate Gowers norm.

This discussion naturally leads to the inverse conjecture \cite{gt,sam, tao-lec}, that is if $(d+1)$-th Gowers norm of $f$ is non-trivial, then $f$ is non-trivially approximable by a degree-$d$ polynomial. This conjecture is easily seen to hold for $d=1$ and has been proved also for $d=2$ \cite{gt,sam}. It is of interest to prove this conjecture for higher values of $d$.

In this paper we show this conjecture, which we will refer to as the 'Inverse Conjecture for the Gowers norm', or, informally, as ICGN, to be false. Let $S_n$ be the elementary symmetric polynomial of degree $n$ in $N$ variables, that is
$$
S_n(x) = \sum_{S \subseteq [N],~|S| = n} \prod_{i\in S} x_i
$$
We prove two claims about symmetric polynomials. Note that here and below a constant is {\it absolute} if it does not depend on $N$.

First, we show Gowers norms of some symmetric polynomials to be non-trivial.
\thm
\label{thm:gowers-positive}
There is an absolute positive constant $\e$ such that for any prime $p$
$$
\|S_{2p}\|_{U^{p+2}} > \e,
$$
Here $S_{2p}$ is viewed as a function over $\FF = \FF_p$.
\ethm

Two versions of this result will be useful later.
\begin{itemize}
\item
A special case $p = 2$.
\beqn
\label{4-positive}
\|S_4\|_{U^4} > \e
\eeqn
\item
An easy generalization: for any $n \ge 2p$,
\beqn
\label{n-p+2-pos}
\|S_n\|_{U^{n-p+2}} > \e
\eeqn
\end{itemize}

In the second claim we show a specific symmetric polynomial to have no non-trivial approximation by polynomials of lower degree.
\thm
\label{thm:non-approx}
Let $p=2$. For any polynomial $g$ of degree $3$ holds
\beqn
\label{4-far}
|\<S_4, g\>| < \exp\{-\alpha N\}
\eeqn
\ethm

We conjecture the second claim of the theorem to be true for any prime number $p$, replacing $3$ with $p+1$ and $4$ with $2p$.

The combination of (\ref{4-positive}) and (\ref{4-far}) shows ICGN to be false for $p = 2$ and $d = 4$.

\subsection{Related work}
Our results have a large overlap with a recent work of Green and Tao \cite{gt07}.

The paper of Green and Tao has two parts. In the first part ICGN is shown to be true when $f$ is itself a polynomial of degree less than $p$ and $d < p$. In the second part, the conjecture is shown to be false in general. In particular the symmetric polynomial $S_4$ is shown to be a counterexample for $p = 2$ and $d = 4$.

To proof of non-approximability of $S_4$ by lower-degree polynomials in \cite{gt07} uses a modification of a Ramsey-type argument due to Alon and Beigel \cite{ab}. Very briefly, this argument shows that if a function over $\FF_2$ has a non-trivial correlation with a multilinear polynomial of degree $d$, then its restriction to a subcube of smaller dimension has a non-trivial correlation with a symmetric polynomial of degree $d$. The problem of inapproximability by symmetric polynomials turns out to be easier to analyze.

This argument gives a somewhat weaker bounds for non-inapproximability of $S_4$, in that it shows $\<S_4,g\> < \log^{-c}(N)$ for any degree-$3$ polynomial $g$ and for an absolute constant $c > 0$.

On the other hand, this argument is more robust than our inapproximability argument. We observe below that it can be readily extended to the case of general prime $p$ and, combined with (\ref{n-p+2-pos}), show ICGN to be false for all $p$.

\subsection{The case of a general prime field}
We briefly observe here that a minor adaptation of the Alon-Beigel argument, together with (\ref{n-p+2-pos}), show the symmetric polynomial $S_{p^2}$ to have a non-negligible $\(p^2\)$-nd Gowers norm over $\FF_p$ and to have no good approximation by lower-degree polynomials. In that, $S_{p^2}$ provides a counterexample to ICGN for any prime $p$.

Indeed, by monotonicity of the Gowers norms (\cite{gt}), and since $p \ge 2$, a direct implication of (\ref{n-p+2-pos}) gives
$$
\|S_{p^2}\|_{U^{p^2}} > \e
$$

On the other hand, let $g$ be a polynomial of degree less than $p^2$ in $N$ variables such that $\<S_{p^2},g\> > \e$. Note that the Alon-Beigel argument (as given in \cite{ab} and in \cite{gt07}) does not seem to be immediately applicable in this case, since $g$ does not have to be multilinear. A way around this obstacle, is to observe, via an averaging argument, that there is a copy of an $N'$-dimensional boolean cube $\{0,1\}^{N'}$, such that restrictions $S'$ and $g'$ of $S_{p^2}$ and of $g$ on this subcube satisfy $\<S', g'\> > \e'$, and $N', \e'$ depend linearly on $N, \e'$. Without loss of generality assume the coordinates of the boolean cube to be $\{1...N'\}$ and consider the functions $S',g'$ as functions in variables $x_1,...,x_{N'}$ (with some fixed assignment of values to variables $x_i,~i > N'$). Now, $S' = \sum_{i=0}^{p^2} a_i S_i$ is a symmetric polynomial of degree $p^2$ over $\FF^{N'}$, with $a_i = 1$, and $g'$ is a polynomial of a degree smaller than $p^2$. Our gain is in that now $g'$ can be replaced by a multilinear polynomial coinciding with $g'$ on the boolean cube, and hence having a non-trivial correlation with $S'$ on the boolean cube.

Now, the Alon-Beigel argument can be applied to show that the symmetric polynomial $S_{p^2}$ has a non-trivial correlation with a symmetric polynomial $h$ of a smaller degree over the boolean cube $\{0,1\}^{N'}$ viewed as a subset of $\FF^{N'}$. This, however, couldn't be true due to a theorem of Lucas, which implies that for a boolean vector $x$ with Hamming weight $w = \sum_{i=1}^{N'} x_i$, the value $S_{p^2}(x)$ depends only on the $3$-rd digit in the representation of $w$ in base $p$, while the value of $h$ depends only on the first $2$ digits.

This completes the argument. We conclude with an observation that this argument directly extends to $S_{p^k}$ for any $k > 1$.

Here is a brief overview of the rest of the paper.
Section~\ref{sec:notions} defines relevant notions and contains proofs of several technical claims.
Theorem~\ref{thm:gowers-positive} is proved in Section~\ref{sec:first-claim}. Theorem~\ref{thm:non-approx} is proved in Section~\ref{sec:second-claim}.

\section{Some useful notions and claims}
\label{sec:notions}

\subsection{Some multilinear polynomials and their properties}
In this sub-section we introduce and discuss certain polynomials over the finite field $\FF$. These polynomials can be conveniently viewed as multi-linear functions on matrices whose entries are elements of $\FF$, or formal variables with values in the field. A basic object we consider is a rectangular $n \times N$ matrix, $N \ge n$. A matrix $M$ with rows $r_1...r_n$ will be denoted by $M[r_1...r_n]$. Sometimes there will be repeated rows. In such a case we consider a partition $\lambda = \(\lambda_1...\lambda_k\)$ of $[n]$, that is $\l_i$ are (possibly empty) subsets of $[n]$, whose disjoint union is $[n]$. We denote by $M_{\lambda}[r_1...r_k]$ the matrix whose rows in positions indexed by elements of $\lambda_i$ equal $r_i$. Note that the partition $\l$ is ordered, in that the ordering of the sets $\l_i$ is relevant. We use the notation $\{\l_1...\l_k\}$ for an unordered partition.

First, we introduce the "symmetric" function $\S$.
We define $\S(M)$ to be the sum of all the permanental  minors of $M$, that is
$$
\S(M) := \sum_{C\subseteq [N],|C| = n} Per\(M_C\),
$$
where $M_C$ is an $n\times n$ submatrix of $M$ which is obtained by deleting all the columns of $M$ except these with indices in $C$.

Let $\l = \(\l_1...\l_k\)$ be a partition of $[n]$, and set $\ell_i = | \l_i |$. Clearly $\S\(M_{\l}\)$ depends only on the cardinalities $\ell_i$ of $\l_i$. This leads to the notation $M\[r_1^{(\ell_1)}...r_k^{(\ell_k)}\]$ which denotes the matrix in which the row $r_1$ appears $\ell_1$ times, followed by $\ell_2$ appearances of the row $r_2$ and so on. In this notation, therefore
$$
\S\(M_{\(\l_1...\l_k\)}[r_1...r_k]\) = \S\(M\[r_1^{(|\l_1|)}...r_k^{(|\l_k|)}\]\)
$$
The second matrix function we consider is the "forward" function $\F$, with
$$
\F(M\[r_1...r_n\])= \sum_{C\subseteq [N],|C| = \{j_1 < j_2 < ... < j_n\}} \prod_{i=1}^n r_i\(j_i\)
$$
Here $r_i(j)$ denote the $j$-th coordinate of the vector $r$.

To connect the two notions, observe that
$$
\S(M[r_1...r_n]) = \sum_{\sigma} \F(M[r_{\sigma_1}...r_{\sigma_n}])
$$
where $\sigma$ runs over all permutations on $n$ items.

The last function we consider is a "hybrid" function $\H$ which has some 'symmetric' and some 'forward' properties. Let $\l = \(\l_1...\l_k\)$ be an ordered partition of $[n]$ with $k$ terms. For another such partition $\theta = \(\theta_1...\theta_k\)$ of $[n]$ write $\theta \sim \l$ if $|\theta_1| = |\l_1|$,...,$|\theta_k| = |\l_k|$.
We define
$$
\H\(M_{\l}[r_1...r_k]\) = \sum_{C\subseteq [N],|C| = \{j_1 < j_2 < ... < j_n\}} \sum_{\theta \sim \l} \prod_{t=1}^k \prod_{i \in \theta_t} r_t\(j_i\)
$$
An alternative view of the functions $\S,\F$ and $\H$ might be helpful at this point. Consider the set of {\it paths} which are one-to-one functions from $[n]$ to $[N]$. Let us call a path $\rho$ monotone on a subset $\{i_1 < i_2 < ... < i_{\ell}\}$ of $[n]$ if $\rho\(i_1\) < \rho\(i_2\) < ... < \rho\(i_{\ell}\)$. A path is (fully) monotone if it is monotone on $[n]$.
Then, for a partition $\l = \(\l_1...\l_k\)$ of $[n]$ and an $n \times N$ matrix $M = M_{\l}$,
$$
\S(M) = \sum_{\mbox{all}~\rho} ~\prod_{i=1}^n M_{i,\rho(i)}
$$
$$
\F(M) = \sum_{\mbox{monotone}~\rho} ~\prod_{i=1}^n M_{i,\rho(i)}
$$
$$
\H(M) = \sum_{\rho~\mbox{monotone on $\l_1...\l_k$}} ~\prod_{i=1}^n M_{i,\rho(i)}
$$
Note that for the function $\H$, similarly to the symmetric function $\S$, holds
$$
\H\(M_{\(\l_1...\l_k\)}[r_1...r_k]\) = \H\(M\[r_1^{(|\l_1|)}...r_k^{(|\l_k|)}\]\)
$$
Observe also that if $\l = \(\{1\}...\{n\}\)$ then $\S(M) = \H(M)$. If $\l = \(\{[n]\}\)$ then $\F(M) = \H(M)$ and $\S(M) = n! \cdot \F(M) = n! \cdot \H(M)$. For a general $\l = \(\l_0...\l_k\)$
\beqn
\label{H-vs-S}
\S(M) = \(\prod_{t=1}^k |\l_t|!\) \cdot \H(M)
\eeqn
Note that this is an identity in $\FF$. In particular, if one of the terms $\l_i$ has cardinality at least $p$ then $\S(M) = 0$ and (\ref{H-vs-S}) provides no information.

To simplify the notation we will usually write $\S(r_1...r_n)$ for $\S(M[r_1...r_n])$, $\F_{\l}(r_1...r_k)$ for $\F\(M_{\l}[r_1...r_k]\)$ and so on.

\subsection{Directional derivatives of symmetric polynomials}
The functions we have defined are relevant to the discussion here for two reasons. First, the elementary symmetric polynomial $S_n(x)$ in $N$ variables can be viewed as the forward function $\F$ applied to the matrix $M[x...x]$, where $M$ has $n$ identical rows equal to $x$. In our notation,
$$
S_n(x) = \F_{\{[n]\}}(x)
$$
Second, it is possible to write a directional derivative $\(S_n\)_{y_1...y_k}$ of $S_n$ of any order as a combination of values of $\F$ on explicitly defined matrices $M$ whose rows are either the indeterminate $x$ or the directions $y_i$.

The basic observation here is the following lemma which is straightforward from the definition of directional derivative.
\lem
\label{lem:basic-der}
Let a polynomial $P(x)$ in $N$ variables be given by
$$
P(x) = \F_{\(\l_0...\l_k\)}(x,y_1...y_k)
$$
Then
$$
P_z(x) = \sum_{A \subset \l_0} \F_{\(A,\l_0 \setminus A,\l_1...\l_k\)}(x,z,y_1...y_k)
$$
In words, when we take the derivative of such a polynomial in direction $z$, we replace some of the rows which contained $x$ with $z$. \elem

As a corollary we have a following expression for higher order derivatives of a symmetric polynomial.
\pro
\label{pro:high-der}
Let $k \le n$, then
$$
\(S_n\)_{y_1...y_k}(x) = \sum_{m=0}^{n-k} ~~\sum_{\ell_1...\ell_k \ge 1,~\sum_i \ell_i = n - m} \H\(x^{(m)},y_1^{(\ell_1)}...y_k^{(\ell_k)}\)
$$
\epro
\prf
Iterating Lemma~\ref{lem:basic-der},
$$
\(S_n\)_{y_1...y_k}(x) = \sum_{\l = \(\l_0,\l_1...\l_k\)} \F_{\l}(x,y_1...y_k)
$$
where the summation is over partitions $\l$ such that $\l_i$ are not empty for $i=1...k$.
Rearranging, this is
$$
\sum_{m=0}^{n-k} ~~~\sum_{\ell_1...\ell_k \ge 1,~\sum_i \ell_i = n - m} ~~~\sum_{\l:~ |\l_0| = m,|\l_1|=\ell_1...|\l_k| = \ell_k}
\F_{\l}(x,y_1...y_k) =
$$
$$
\sum_{m=0}^{n-k} ~~\sum_{\ell_1...\ell_k \ge 1,~\sum_i \ell_i = n - m} \H\(x^{(m)},y_1^{(\ell_1)}...y_k^{(\ell_k)}\)
$$
\eprf

We can give explicit expressions for the coefficients of $\(S_n\)_{y_1...y_k}(x)$. Fix $m$ indices $j_1 < j_2 < ... < j_m$ for $0 \le m \le n-k$, and let $a$ be the coefficient of $x_{j_1} \cdots x_{j_m}$ in $\(S_n\)_{y_1...y_k}$.
\cor
\label{cor:coef-high-der}
\begin{itemize}
\item
$$
a = \sum_{\ell_1...\ell_k \ge 1,~\sum_i \ell_i = n - m} \H^{\{j_1...j_m\}}\(y_1^{(\ell_1)}...y_k^{(\ell_k)}\)
$$
\item
If $k + m + p > n+1$ then
$$
a = \sum_{\ell_1...\ell_k \ge 1,~\sum_i \ell_i = n - m} \(\prod_{i=1}^k \ell_i!\)^{-1} \cdot \S^{\{j_1...j_m\}}\(y_1^{(\ell_1)}...y_k^{(\ell_k)}\)
$$
\end{itemize}
Here, for a subset of indices $T \subseteq [N]$,  $\H^T(M)$ returns the value of the matrix function $\H$ applied to the $n \times (N - |T|)$ matrix obtained from $M$ by deleting columns in $T$. The function $\S^T(M)$ is defined similarly.
\ecor
\prf
The first claim is immediate from Proposition~\ref{pro:high-der}.
The second claim follows from the first claim, from (\ref{H-vs-S}), and from the simple observation that if $k + m + p > n + 1$ then $\ell_i < p$ for $i = 1...k$ in the above summation, which means $\ell_i !$ is invertible in $\FF_p$.
\eprf

\xam
\label{xam:4-2}
The following "toy" example will be relevant for the case of the binary field. It is sufficiently simple to illustrate what's going on behind the cumbersome formulas. Consider $P = \(S_4\)_{y,z}$. Then $P$ is a quadratic polynomial and for $1 \le i < j \le N$
$$
\mbox{coef}_{x(i) x(j)}(P) = \sum_{k \not = l,~k,l \not \in \{i,j\}} y(k) z(l) = \S^{\{i,j\}}\(y,z\)
$$
\exam
Continuing with the same example, note that it convenient to express the symmetric function $\S\(y,z\)$ via inner products of vectors $y,z,\1$, where $\1$ is the all-$1$ vector of length $N$.
$$
\S(y,z) = \sum_{k \not = l} y(k) z(l) = \<y,\1\> \cdot \<z,\1\> - \<yz,\1\>
$$
Here we take $yz$ to be the vector whose coordinates are point-wise inner products of the coordinates of $y$ and $z$, that is $(yz)(i) = y(i) z(i)$. Of course, $\<yz,\1\>$ is the same as $\<y,z\>$.

Similarly, we can express the 'incomplete' symmetric function $\S^{\{i,j\}}(y,z)$ via the complete symmetric function $\S(y,z)$
minus forbidden terms, as follows
$$
\S^{\{i,j\}}(y,z) =  \S(y,z) - \Big (z(i) + z(j) \Big ) \<y,\1\> - \Big (y(i) + y(j)\Big ) \<z,\1\> + \Big (y(i) z(j) + y(j) z(i)\Big )
$$
Note the "inclusion-exclusion" structure in the two expressions above. (To make it even clearer we use "+" and "-" notation, though in the binary field both are, of course, the same.) This structure becomes more evident as we pass to our next order of business, which is expressing, for general $n$ and $k$, the coefficients of $\(S_n\)_{y_1...y_k}$ via inner products of vectors $y_1...y_k,\1$.

\subsection{Inclusion-Exclusion formulas for symmetric functions}
Some notation: Given $m$ vectors $y_1...y_m$ and a subset $\T \subseteq [m]$, let $y_{\T}$ to be vector whose coordinates are point-wise products of the corresponding coordinates of $y_i$, $i \in \T$. Let $\S\(y[\T]\)$ for the value of the function $\S$ on a matrix with $|\T|$ rows $y_i$, $i\in \T$. Let $\<y_{\T}\>$ be the polynomial $\<y_{\T},\1\> = \sum_{j=1}^N \prod_{i \in \T} y_i(j)$.

We start with an auxiliary lemma expressing the incomplete symmetric function $\S^{\{k\}}\(r_1...r_n\)$ as a polynomial in the $k$-th coordinate of the vectors $r_i$ and in complete symmetric functions applied to sub-matrices of $M[r_1...r_n]$.
\lem
\label{lem:incomplete-sym}
$$
\S^{\{k\}}\(r_1...r_n\) = \sum_{\T \subseteq [n]} (-1)^{|\T|} (|\T|)! \cdot r_{\T}(k) \cdot \S\(r\Big[[n] \setminus \T\Big]\)
$$
From now on we assume $r_{\emptyset}$ to be the all-$1$ vector, and $\S\(r[\emptyset]\)$ to equal $1$.
\elem
\prf
The proof is by induction on $n$. For $n=1$ both sides equal $\sum_{j=1}^N r_1(j) - r_1(k)$.

For $n > 1$, observe that
$$
\S^{\{k\}}\(r_1...r_n\) = \S\(r_1...r_n\) - \sum_{i=1}^n r_i(k) \cdot \S^{\{k\}}\(r\Big[[n]\setminus \{i\}\Big]\)
$$
and the claim is easily verified using the induction hypothesis.
\eprf

Now we can state two main claims of this section. The first expresses the complete symmetric function $\S\(r_1...r_n\)$ via inner products $\<r_T\>$.
\pro
\label{pro:sym}
$$
\S\(r_1...r_n\) = \sum_{\l = \{\l_1...\l_m\}} \prod_{t=1}^m \((-1)^{|\l_t|-1} \(|\l_t| - 1\)! \cdot \<r_{\l_t}\>\)
$$
In this summation $\l = \{\l_1...\l_m\}$ runs over all unordered partitions of $[n]$ with non-empty $\l_i$.
\epro
\prf
Again, the proof is by induction on $n$. For $n=1$ both sides equal $\sum_{j=1}^N r_1(j)$. For $n > 1$ we have
$$
\S(r_1...r_n) = \sum_{k=1}^N r_n(k) \cdot \S^{\{k\}}(r_1...r_{n-1})
$$
Using Lemma~\ref{lem:incomplete-sym} and the induction hypothesis,
$$
\S(r_1...r_n) = \sum_{k=1}^N r_n(k) \cdot \sum_{\T\subseteq [n-1]} (-1)^{|\T|} (|\T|)! \cdot r_{\T}(k) \cdot \S\(r\Big[[n-1]\setminus \T\Big]\) =
$$
$$
\sum_{\T\subseteq [n-1]} (-1)^{|\T|} (|\T|)! \cdot \<r_{\T \cup [n]}\> \cdot \S\(r\Big[[n-1]\setminus \T\Big]\)
$$
Consider the summand corresponding to $\T = [n-1]$. Recall the boundary assumption $\S\(r[\emptyset]\) = 1$. Hence this summand is $(-1)^{n-1} (n-1)! \cdot \<r_{[n]}\>$. This summand therefore corresponds to the partition $\l = \{[n]\}$ in the claim of the proposition.

For $\T$ a proper subset of $[n-1]$, we use the induction hypothesis to obtain
$$
\S(r_1...r_n) =  \sum_{\T\subseteq [n-1]} (-1)^{|\T|} (|\T|)! \cdot \<r_{\T \cup [n]}\> \cdot
\sum_{\theta = \{\theta_1...\theta_l\}}  \prod_{t=1}^l \((-1)^{|\theta_t| - 1} \(|\theta_t| - 1\)! \cdot \<r_{\theta_t}\>\) +
$$
$$
(-1)^{n-1} (n-1)! \cdot \<r_{[n]}\>
$$
Here $\theta$ runs over all the unordered partitions of $[n-1] \setminus \T$ with non-empty $\theta_i$. Observe that each pair $(\T,\theta)$ leads to a unique partition $\l = \{\l_1...\l_{l+1}\} = \{\theta_1...\theta_l,\T \cup [n]\}$ of $[n]$. Rearranging the terms, the last summation can be written as
$$
\sum_{\l = \(\l_1...\l_m\)} \prod_{t=1}^m \((-1)^{|\l_t|-1} \(|\l_t| - 1\)! \cdot \<r_{\l_t}\>\)
$$
completing the proof of the proposition.
\eprf

The second claim expresses the incomplete symmetric function $\S^{\{j_1...j_k\}}\(r_1...r_n\)$ as a polynomial in the missing coordinates $j_1...j_k$  of the vectors $r_i$ and in complete symmetric functions applied to sub-matrices of $M[r_1...r_n]$.
Note that Lemma~\ref{lem:incomplete-sym} is a special case $k=1$ of this claim.

\pro
\label{pro:incomp-sym}
$$
\S^{\{j_1...j_k\}}\(r_1...r_n\) = \sum_{\T = \(\T_1...\T_k\)} \prod_{t = 1}^k \((-1)^{|\T_t|} (|\T_t|)! \cdot r_{\T_t}\(j_t\)\) \cdot \S\(r\Big[[n] \setminus \cup_t \T_t\Big]\)
$$
Here the summation is on all ordered set systems $\tau$ such that the terms $\tau_t$ are disjoint subsets of $[n]$. The terms may also be empty.
\epro
\prf
The proof is by induction on $k$ and $n$. The case $k = 1$ is treated in Lemma~\ref{lem:incomplete-sym}.

Consider the case $n = 1$. On one hand $\S^{\{j_1...j_k\}}\(r_1\) = \sum_{j=1}^N r_1(j) - \sum_{t=1}^k r_1\(j_t\)$.
We claim that this value can be also represented as
$$
\sum_{\T = \(\T_1...\T_k\)} \prod_{t = 1}^k \((-1)^{|\T_t|} (|\T_t|)! \cdot r_{\T_t}\(j_t\)\) \cdot \S\(r\Big[[1] \setminus \cup_t \T_t\Big]\)
$$
Here $\T_i$ are disjoint subsets of $[1]$. Observe that there are $k+1$ summands in this expression, corresponding to different set systems $\T$. Let $\T^{(0)}$ denote the set system with $k$ empty terms, and let $\T^{(t)}$, for $t = 1...k$ denote the set system with $\T_t = \{1\}$ and all the remaining terms are empty. The summand corresponding to $\T^{(0)}$ is $\S\(r_1\) = \sum_{j=1}^N r_1(j)$. The summand corresponding to $\T^{(t)}$ is $\(- r_1\(j_t\)\) \cdot \S\(r_{\emptyset}\) = - r_1\(j_t\)$, and we are done in this case.

For $k, n > 1$, we have
$$
\S^{\{j_1...j_k\}}\(r_1...r_n\) = \S^{\{j_1...j_{k-1}\}}\(r_1...r_n\) - \sum_{i=1}^n r_i\(j_k\) \cdot \S^{\{j_1...j_k\}}\(r\Big[[n]\setminus \{i\}\Big]\)
$$
By the induction hypothesis, this is
$$
\sum_{\theta = \(\theta_1...\theta_{k-1}\)} \prod_{t = 1}^{k-1} \((-1)^{|\theta_t|} (|\theta_t|)! \cdot r_{\theta_t}\(j_t\)\) \cdot \S\(r\Big[[n] \setminus \cup_t \theta_t\Big]\) -
$$
$$
\sum_{i=1}^n r_i\(j_k\) \cdot \sum_{\mu^{(i)} = \(\mu^{(i)}_1...\mu^{(i)}_k\)} \prod_{u = 1}^k \((-1)^{|\mu^{(i)}_u|} (|\mu^{(i)}_u|)! \cdot r_{\mu^{(i)}_u}\(j_u\)\) \cdot \S\(r\Big[[n] \setminus \cup_t \mu^{(i)}_t \setminus \{i\}\Big]\)
$$
Here the summation is on all ordered set systems $\theta$ such that the terms $\theta_t$ are disjoint subsets of $[n]$ and on
ordered set systems $\mu^{(i)}$, $i=1...n$ such that the terms $\mu^{(i)}_u$ are disjoint subsets of $[n] \setminus \{i\}$.

Given a set system $\theta = \(\theta_1...\theta_{k-1}\)$ we define a set system $\T = \(\T_1...\T_k\)$ by setting $\T_t = \theta_t$, $t = 1...k-1$ and $\T_k = \emptyset$. Given a set system $\mu^{(i)} = \(\mu^{(i)}_1...\mu^{(i)}_k\)$ we define a set system $\T = \(T_1...\T_k\)$ by setting $\T_u = \mu^{(i)}_u$, $u = 1...k-1$ and $\T_k = \mu^{(i)}_k \cup \{i\}$. In both cases we have obtained a set system of the type we want, that is an ordered family of $k$ disjoint subsets of $[n]$. Moreover, each such system with empty $k$-th term is obtained exactly once, from the corresponding $\theta$-system, and each system with non-empty $k$-th term $\T_k$ is obtained exactly $|\T_k|$ times, from systems $\mu^{(i)}$ with $i \in \T_k$. Rearranging the terms and the signs, the last expression is precisely
$$
\sum_{\T = \(\T_1...\T_k\)} \prod_{t = 1}^k \((-1)^{|\T_t|} (|\T_t|)! \cdot r_{\T_t}\(j_t\)\) \cdot \S\(r\Big[[n] \setminus \cup_t \T_t\Big]\),
$$
completing the proof.
\eprf

\subsection{Some properties of Gowers' norms}
The main result in this subsection shows that if a function from $\FF^N$ to $\FF$ is fixed on a subset of $\FF^N$ defined by low-degree polynomial constraints, then it has a non-trivial Gowers norm of an appropriate order.

Recall that for a vector $x \in \FN$, $x^i$ stands for a vector in $\FN$ whose coordinates are $i$-th powers of the coordinates of $x$.

\pro
\label{pro:Gow-high}
Let $K$ be an absolute constant. Let $y_{i,j}$, $i = 1...p-1$, $j = 1...K$, be $K(p-1)$ vectors in $\FN$. Let $M$ be a subset of $\FN$ defined by the constraints $\<x^i,y_{i,j}\> = 0$ for all $i,j$.

Let $f$ be a function from $\FN$ to $\FF$. Assume that $f$ is fixed on $M$. Then
$$
\|f\|_{U^p} > \(\frac{|M|}{2^N}\)^2 =: Pr^2\{M\}
$$
\epro
\prf
Let $f_{|M} \equiv c_0$.

Consider a subspace $V$ of polynomials of degree at most $p-1$ in $\FF[x_1...x_N]$ spanned by the polynomials $\<x^i,y_{i,j}\>$, for all $i,j$. We will first find a polynomial $g \in V$ such that $|\<f,g\>| \ge Pr\{M\}$. This, combined with a lemma from \cite{gt}, will imply the claim of the proposition.

Let $\b = \(b_{i,j}\)$, $i = 1...p-1$, $j = 1...K$, be a matrix with entries in $\FF$. Let $c \in \FF$. Set
$$
\mu(\b,c) = Pr\Big\{x:~f(x) = c ~\wedge~ \<x^i,y_{i,j}\> = b_{i,j}~\mbox{for all}~i,j\Big\}
$$
Note that, by assumption, for a zero matrix $\b$ holds $\mu\(\b,c_0\) = Pr\{M\}$. In other words, $\mu(\b,c) = 0$ and for $\b  = 0$ any $c \not = c_0$.

Now, for any $g(x) = \sum_{i,j} a_{i,j} \<x^i,y_{i,j}\>$ in $V$ holds
$$
\<f,g\> = \E e(f-g) = \sum_{\b,c} \mu(\b,c) \cdot e\( c - \<\a,\b\>\)
$$
where $\a = \(a_{i,j}\)_{i,j}$ and $\<\a,\b\> = \sum_{i,j} a_{i,j} b_{i,j}$. Averaging over $V$, we have
$$
\E_{g \in V} \<f,g\> = \frac{1}{|V|} \sum_{\a} \sum_{\b,c} \mu(\b,c) \cdot e\(c - \<\a,\b\>\) =
\frac{1}{|V|} \sum_{\b,c} \mu(\b,c) \cdot e(c) \sum_{\a} e\(-\<\a,\b\>\) =
$$
$$
 \sum_c \mu(0,c) \cdot e(c) = \mu\(0,c_0\) \cdot e\(c_0\) = Pr\{M\} \cdot e\(c_0\)
$$
This means, there is $g \in V$ with $|\<f,g\>| \ge Pr\{M\}$. We conclude the proof of the proposition by quoting a lemma from \cite{gt}, which states that $|\<f,g\>| \ge \e$ implies $\|f\|_{U^p} \ge \e$.
\eprf

\subsection{Asymptotic uniformity and independence of some random variables}
In this subsection we deal with another property of multiviarite polynomials. Let $n$ be fixed integer and let $N$ be an integer parameter growing to infinity. Let $r_1...r_n$ be $n$ vectors in $\FF^N$. Let $\kappa = \(k_1...k_n\)$ be a non-zero sequence of integers $0 \le k_i < p$. For each such sequence define a polynomial $X_{\kappa}\(r_1,...,r_n\) = \sum_{j=1}^N \prod_{i=1}^n r^{k_i}_i(j)$.

Now,let $r_1...r_n$ be chosen uniformly and independently from $\FF^N$. We claim that for a large $N$ the random variables $X_{\kappa}\(r_1,...,r_n\)$ are nearly independent and uniformly distributed over $\FF$. Let  $X = \(X_{\kappa}\)_{\kappa}$, and let
$K = p^n$.
\pro
\label{pro:asymp-ind}
Let $U$ be the uniform distribution on $\FF^K$. Let $P$ be distribution of $X$ on $\FF^K$.
Let $\|\cdot \|$ denote the statistical ($l_1$) distance between distributions.

Then there is a constant $c > 0$ depending on $n,p$ but not on $N$ such that
$$
\|P - U\| \le \exp\left\{-cN\right\}
$$
\epro
\prf
We start from a simple observation that Fourier transform of a uniform distribution is the delta function at $0$. In addition, the two following statements are equivalent up to constants: 'a distribution is exponentially close to uniform' and 'all non-zero Fourier coefficients of the distribution are exponentially close to zero'. Accordingly, we will show that all the non-zero Fourier coefficients of $P$ tend exponentially fast in $N$ to zero.

Consider a character $\chi(y) = \xi^{\<y,a\>}$, corresponding to a non-zero vector $a = \(a_{\kappa}\)_{\kappa} \in \FF^K$. (Recall that $\xi = e^{2\pi i/p}$ is the $p$-th primitive root of unity.) Then, normalizing appropriately,
$$
\widehat{P}(\chi) = \sum_y P(y) \bar{\chi}(y) = \sum_y Pr\{X = y\} \cdot \xi^{-\sum_{\kappa} a_{\kappa} y_{\kappa}} = \E \xi^{-\sum_{\kappa} a_{\kappa} X_{\kappa}}
$$
Let $P_a$ denote the distribution of the random variable $X_a = \sum_{\kappa} a_{\kappa} X_{\kappa}$. Then we have shown $\widehat{P}(\chi) = \widehat{P_a}(1)$. We will show the non-zero Fourier coefficients of $P_a$ to be exponentially small, completing the proof of the proposition.

We have
$$
X_a\(r_1,...,r_n\) =  \sum_{\kappa} a_{\kappa} P_{\kappa}\(r_1,...,r_n\) = \sum_{j=1}^N ~~\sum_{\kappa = \(k_1...k_n\)} a_{\kappa}~~\prod_{i=1}^n r^{k_i}_i(j)
$$
Let $x_i$ be elements of the field $\FF$. Consider an $n$-variate polynomial
$$
Q(x_1...x_n) = \sum_{\kappa = \(k_1...k_n\)} a_{\kappa} ~~\prod_{i=1}^n x_i^{k_i}
$$
Since not all of the coefficients $a_{\kappa}$ are zero, and since all $\kappa$ are non-zero sequences, $Q$ is a multi-variate polynomial of degree at least $1$ in $\FF[x_1...x_n]$, and therefore attains at least two values with probability bounded away from zero. Now, $X_a = \sum_{j=1}^N Q\(r_1(j)...r_n(j)\)$ is a sum of $N$ independent copies of $Q$. Let $\mu$ denote the distribution of $Q$ on $\FF$. Then the distribution $P_a$ of $X_a$ is $\mu^{\ast N}$, the $N$-wise convolution of $\mu$ with itself. Since $p$ is prime, $\widehat{\mu}(0) = 1$, and $|\widehat{\mu}| < 1$ everywhere else. Therefore,
$\widehat{P_a} = \(\widehat{\mu}\)^N$ tends to the delta function at $0$ exponentially fast in $N$, completing the proof.
\eprf

\subsection{Estimates on the number of common zeroes of some families of polynomials}
The main claim of this subsection is the following proposition.
\pro
\label{pro:commonz}
Let $M$ be the ring of $\FF$-valued functions on $\FN$, that is $M = \FF[x_1...x_N]/I$, where $I$ is the ideal $\(x^p_1 - x,...,x^p_N -x\)$. Let $f_1...f_K$ be polynomials in $M$. Let $S$ be the set of common zeroes of $f_1...f_K$, that is
$$
S = \Big\{u \in \FN:~f_1(u)=...=f_K(u) = 0 \Big\}
$$
Then
$$
|S| \le \mbox{dim} \(M/J\)
$$
where $J$ is the ideal generated by $\{f_i\}$, and $\mbox{dim} \(M/J\)$ denotes the dimension of $\mbox{dim} \(M/J\)$, viewed as a vector space over $\FF$.
\epro
\prf
For each $u \in S$, let $q_u \in M$ be defined by $q_u(u) = 1$ and $q_u(v) = 0$ for all $v \not = u$. We will show that the family $\{q_u + J\}_{u \in S}$ is linearly independent in $M/J$. This will immediately imply the claim of the proposition.

Consider a linear combination $q = \sum_{u \in S} \lambda_u q_u$ such that $q \in J$. Let $v \in S$. We compute $q(u)$ in two ways. First, since $q \in J$, we have $q(v) = 0$. On the other hand, $q(v) = \sum_{u \in S} \lambda_u q_u(v) = \lambda_v$. This shows $\lambda_v = 0$ for all $v \in S$, completing the proof.
\eprf

In some cases, the dimension of $M/J$ is easy to estimate.
\lem
\label{lem:example}
Let $p = 2$, let $K = {N \choose k}$, and let $\{f_I\}$ be indexed by $k$-subsets $I$ of $[N]$. Assume that for any such subset $I$ holds
\beqn
\label{cond:lower-deg}
\mbox{deg}\(f_I(x) - \prod_{i\in I} x_i\) \le k - 1
\eeqn
Then,
$$
\mbox{dim} \(M/J\) \le \sum_{j=0}^{k-1} {N \choose j}
$$
\elem
\prf
We will construct a generating subset of the vector space $M/J$ of cardinality at most $\sum_{j=0}^{k-1} {N \choose j}$. We start from a trivial generating set $\{m + J\}$, where $m$ runs through all the $2^N$ multi-linear monomials in $N$ variables. Now, in the factor space $M/J$, we can replace any product of $k$ variables, $\prod_{i\in I} x_i$, by a polynomial of degree smaller than $k$. Iterating this procedure, we arrive to a generating set spanned by $\{s + J\}$, where $s$ now runs through  $\sum_{j=0}^{k-1} {N \choose j}$ monomials of degree at most $k-1$.
\eprf

\section{Proof of Theorem~\ref{thm:gowers-positive}}
\label{sec:first-claim}
We need to show that
$$
\|S_{2p}\|_{U^{p+2}} > \e
$$
for an absolute constant $\e$.

We remark that (\ref{n-p+2-pos}) can be shown exactly in the same way, replacing $2p$ with $n$ and $p + 2$ with $n - p + 2$ throughout.

Recall (\cite{gt}) that $\|f\|_{U^{p+2}} = \E^{1/2^{p+2}}_{y,z} \|f_{y,z}\|^{2^p}_{U^p}$. Since the Gowers' norms are nonnegative, it will suffice to show that $\|f_{y,z}\|_{U^p}$ is non-negligible for a non-negligible fraction of directions $y,z$.

Let
$$
A = \Big\{(y,z):~\<y^a,z^b\> = 0~\mbox{for all}~0 \le a,b < p\Big\}
$$
By Proposition~\ref{pro:asymp-ind}, for uniformly and independently chosen directions $y,z$, and for a sufficiently large $N$, the probability of $A$  is very close to $p^{-p^2}$. Therefore, $A$ is a non-negligible event. We will now show that for any $(y,z) \in A$ holds $\|f_{y,z}\|_{U^p} > \e'(y,z)$, for an appropriate function $\e'$.

Fix $(y,z)$ in $A$. Let $f = \(S_{2p}\)_{y,z}$. Let
$$
M = M(y,z) = \Big\{x:~\<x^i,y^a z^b\> = 0~\mbox{for all}~1 \le i\le p-1,0 \le a,b < p\Big\}
$$
We will show that $f$ is fixed on $M$. Assuming this, by Proposition~\ref{pro:Gow-high}, we have $\|f_{y,z}\|_{U^p} > Pr^2\{M\}$, and therefore
$$
\|f\|^{2^{p+2}}_{U^{p+2}} = \E_{y,z} \|f_{y,z}\|^{2^p}_{U^p} \ge Pr\{A\} \cdot \E_{(y,z) \in A} Pr^{2^{p+1}}\{M(y,z)\} \ge
$$
$$
Pr\{A\} \cdot \E^{2^{p+1}}_{(y,z) \in A} Pr\{M(y,z)\} \ge \(Pr\{A\} \cdot \E_{(y,z) \in A} Pr\{M(y,z)\}\)^{2^{p+1}} =
$$
$$
Pr^{2^{p+1}}\Big\{x:~\<x^i y^a z^b\> = 0~\mbox{for all}~0 \le a,b,i \le p-1\Big\} \ge \Omega\(p^{-p^3 \cdot 2^{p+1}}\)
$$
The last inequality follows from Proposition~\ref{pro:asymp-ind}, since random variables $\<x^i y^a z^b\>$ are asymptotically uniform and independent.

It remains to prove the following fact.
\lem
\label{lem:vanish}
Let $x,y,z$ be three vectors in $\FN$ satisfying $\<x^i y^a z^b\> = 0~\mbox{for all}~0 \le a,b,i \le p-1$. Then
$$
\(S_{2p}\)_{y,z}(x) = \H\(y^{(p)},z^{(p)}\)
$$
\elem
\prf
By Proposition~\ref{pro:high-der},
$$
\(S_{2p}\)_{y,z}(x) = \sum_{m=0}^{2p-2} ~~\sum_{a,b \ge 1,~a + b = 2p - m} \H\(x^{(m)},y^{(a)},z^{(b)}\)
$$
We claim that all of the summands on the right, except (possibly) $\H\(y^{(p)},z^{(p)}\)$ are $0$.

There are two possible cases to consider. The easier case is when $a,b,m < p$. In such a case, by (\ref{H-vs-S}), $\H\(x^{(m)},y^{(a)},z^{(b)}\)$ is proportional to $\S\(x^{(m)},y^{(a)},z^{(b)}\)$.
By Proposition~\ref{pro:sym}, the symmetric function $\S\(x^{(m)},y^{(a)},z^{(b)}\)$ is a polynomial in $\<x^i y^a z^b\>$, which vanishes when all of these inner products are $0$.

In the second case, one of the indices $a,b,m$ is at least $p$. Note, that there could be at most one such index (barring the case $a = b = p$). We may assume this index is $m$. We claim that in this case $\H\(x^{(m)},y^{(a)},z^{(b)}\)$ can be written as a linear combination of hybrid functions $\H\(x^{(\ell)},r_1,...,r_{m-\ell}\)$, where $\ell < m$ and the vectors $r_i$ are of the form $x^{\alpha} y^{\beta} z^{\gamma}$. Note that this will suffice to prove the lemma, since iterating this step will express $\H\(x^{(m)},y^{(a)},z^{(b)}\)$ as a linear combination of symmetric functions in $r_i$, and these functions vanish.

Consider $\H\(x^{(m)},y^{(a)},z^{(b)}\)$. For notational convenience, let $w_1...w_{a+b}$ stand for the vectors $y...y,z...z$ ($y$ taken $a$ times and $z$ taken $b$ times). Note that both $a$ and $b$ are smaller than $p$. Using Corollary~\ref{cor:coef-high-der} and Proposition~\ref{pro:incomp-sym},
$$
\H\(x^{(m)},y^{(a)},z^{(b)}\) = \(a! \cdot b!\)^{-1} \cdot \sum_{i_1<i_2<...<i_m} x_{i_1} x_{i_2}\cdots x_{i_m} \S^{\{i_1...i_m\}} \(y^{(a)},z^{(b)}\) =
$$
$$
\(a! \cdot b!\)^{-1} \cdot \sum_{i_1<i_2<...<i_m} x_{i_1} x_{i_2}\cdots x_{i_m} \cdot \sum_{\T = \(\T_1...\T_m\)} ~\prod_{t = 1}^m \((-1)^{|\T_t|} (|\T_t|)! \cdot w_{\T_t}\(i_t\)\) \cdot \S\(w\Big[[a+b] \setminus \cup_t \T_t\Big]\)
$$
Here the inner summation is on all ordered set systems $\tau$ such that the terms $\tau_t$ are disjoint subsets of $[a+b]$. The terms may also be empty.

Let us attempt to simplify the double summation we obtained. First, we may disregard the constant term $\(a! \cdot b!\)^{-1}$. Next,
observe that, as before, all symmetric functions of the form $\S\(w[T]\)$ vanish, unless $T$ is empty, in which case they equal $1$. Therefore, we may consider the double summation
$$
\sum_{i_1<i_2<...<i_m} x_{i_1} x_{i_2}\cdots x_{i_m} \cdot \sum_{\T = \(\T_1...\T_m\)} ~\prod_{t = 1}^m \((-1)^{|\T_t|} (|\T_t|)! \cdot w_{\T_t}\(i_t\)\)
$$
Here the inner summation is on all ordered partitions $\tau$ of $[a+b]$. The terms $\tau_t$ may also be empty. Changing the order of summation, and ignoring the constant term $(-1)^{a+b}$, we get
$$
\sum_{\T = \(\T_1...\T_m\)} ~\prod_{t = 1}^m  (|\T_t|)! \cdot  \sum_{i_1<i_2<...<i_m} ~\prod_{t=1}^m  \(x\cdot w_{\T_t}\)\(i_t\) =
\sum_{\T = \(\T_1...\T_m\)} ~\(\prod_{t = 1}^m  (|\T_t|)!\) \cdot \F\(xw_{\T_1},xw_{\T_2},...,xw_{\T_m}\)
$$
Consider the last expression. Let us use some more notation. For an ordered partition $\T = (\T_1...\T_m)$, let $n = n(\T)$ be the number of empty terms. Let $\{\T_1...\T_m\}$ denote the unordered version of this partition, where the first $n(\T)$ terms are taken, by agreement, to be the empty ones. Then we can rewrite this expression as
$$
\sum_{\T = \{\T_1...\T_m\}} ~\(\prod_{t = 1}^m  (|\T_t|)!\) \cdot \H\(x^{(n)},xw_{\T_{n+1}},...,xw_{\T_m}\)
$$
Now, clearly not all the terms in the partition are empty and, therefore, $n(\T) < m$ for all $\T$, completing the proof of our last claim, of the lemma, and of the theorem.
\eprf

\section{Proof of Theorem~\ref{thm:non-approx}}
\label{sec:second-claim}

Let $p=2$. We will show there is an absolute constant $\alpha > 0$ such that for any polynomial $g$ of degree at most $3$ in $N$ variables holds
$$
\<S_4, g\> <\exp \{-\alpha N\}
$$

A first step is to observe that there is a relation between the inner product of two functions and the average inner product of their derivatives.
\lem
\label{lem:mixed-der}
For any two functions $f$ and $g$ holds
$$
\<f,g\>^4  \le \E_y \<f_y,g_y\>^2
$$
\elem
\prf
This is an immediate corollary of a lemma in \cite{sam}, but we give the elementary proof for completeness. By the Cauchy-Schwarz inequality,
$$
\E_y \<f_y,g_y\>^2 \ge \E^2_y \<f_y,g_y\> = \E^2_{x,y} (-1)^{f(x) + f(x+y) + g(x) + g(x+y)} = \E^4 (-1)^{f(x) + g(x)} = \<f,g\>^4
$$
\eprf
\cor
$$
\<f,g\>^8 \le \E_{y,z} \<f_{y,z},g_{y,z}\>^2
$$
\ecor
We will show that for any polynomial $g$ of degree at most $3$ holds $\E_{y,z} \<\(S_4\)_{y,z},g_{y,z}\>^2 \le \exp\{-\alpha N\}$. First, here is a brief overview of the argument.

The point is that taking second derivatives makes life easier, since a second derivative of $g$ is a linear function, and a second derivative of $S_4$ is a quadratic. We therefore need to show that for the large majority of directions $y,z$, the quadratic function $\(S_4\)_{y,z}$ has a small inner product with the linear function $(-1)^{g_{y,z}}$. In this we will be helped by a theorem of Dixon giving a structural description of quadratic polynomials, which, in particular, characterizes the Fourier transform  of functions of the type $(-1)^Q$, where $Q$ is a quadratic. In fact, setting $Q = \(S_4\)_{y,z}$ we will see that for many of the directions $y,z$ the Fourier coefficients of $(-1)^Q$ will be exponentially small. For the remaining directions, these Fourier coefficients
will be supported on an explicit easy to describe $3$-dimensional affine subspace depending on $y,z$. We will then argue that for any fixed polynomial $g$ of lower degree, the support of the character $(-1)^{g_{y,z}}$ lies in this affine subspace with exponentially small probability over $y,z$.

We proceed with computing the second derivative $Q = \(S_4\)_{y,z}$.

\subsection{Second derivatives of $S_4$}
Write
$
Q(x) = \sum_{i < j} q_{i,j} x(i) x(j) + \sum_i  \ell_i x(i) + c.
$

By Proposition~\ref{pro:high-der} or by Example~\ref{xam:4-2}.
$$
q_{i,j} =  \S(y,z) - \<y,\1\> \cdot \Big(z(i) + z(j)\Big) + \<z,\1\> \cdot \Big(y(i) + y(j)\Big)  + \Big(y(i) z(j) + y(j) z(i)\Big)
$$
At this point we invoke (a corollary of) a theorem of Dixon \cite{MS}:
\thm
\label{Dixon}
Let $Q(x) = \sum_{i < j} q_{i,j} x(i) x(j) + \sum_i  \ell_i x(i) + c$ be a quadratic polynomial over $\FF_2$. Consider the symmetric matrix with zeros on the diagonal and off-diagonal entries given by $S_{i,j} = S_{j,i} = q_{i,j}$. Let the rank of $B = 2h$ (it is always even). Then the function $(-1)^Q$ has $2^{2h}$ non-zero Fourier coefficients of absolute value $2^{-h}$. Moreover, all these coefficients lie in an $2h$-dimensional affine subspace of $\FF_2^n$.
\ethm

Consider the matrix $B$ in our case. Some notation: let $J$ be the matrix with $0$ on the diagonal and $1$ off the diagonal. Let $u \otimes v$ denote the outer product $u v^t$. Then,
$$
B = \S(y,z) \cdot J + \<y,\1\> \cdot \Big(z\otimes \1 + \1 \otimes z\Big) + \<z,\1\> \cdot \Big(y\otimes \1 + \1 \otimes y\Big) + \Big(y \otimes z + z \otimes y\Big)
$$
Since the rank of $J$ is at least $N-1$ and the rank of the remaining matrices is at most $2$, the matrix $B$ is almost of full rank if $\S(y,z) = 1$. In this case, by Theorem~\ref{Dixon}, the Fourier coefficients of $(-1)^Q$ are exponentially small.

We therefore may assume $\S(y,z) = 0$. In this case the quadratic part of $Q$ may be written as
$$
\sum_{i < j} q_{i,j} x(i) x(j) = \<y,\1\> \cdot \<x,\1\>\<x,z\> + \<z,\1\> \cdot \<x,\1\> \<x,y\> + \Big(\<x,y\> \<x,z\> + \<x,yz\>\Big)
$$
Recall that $yz$ denotes the pointwise product of vectors $y$ and $z$.

This implies the non-zero Fourier coefficients of $\sum_{i < j} q_{i,j} x(i) x(j)$ lie in a $3$-dimensional affine subspace of $\FF_2^n$. The linear part of this subspace is spanned by the vectors $y,z,\1$ and it is shifted by a vector $yz$.

Next, consider the linear part $\sum_i \ell(i) x(i)$ of $Q$. By Proposition~\ref{pro:high-der},
$$
\ell(i) = \H^{\{i\}}\(y^{(2)},z\) + \H^{\{i\}}\(y,z^{(2)}\) =
$$
$$
\sum_{j < k < l \not = i} \Big(y(k) y(l) z(j) + y(j) y(l) z(k) +
y(j) y(k) z(l)\Big) +
\Big(y(j) z(k) z(l) + y(k) z(j)  z(l) + y(l) z(j) z(k)\Big)
$$
This can be directly verified to be equal to
$$
\Big(\S(y,z) + \S(z,z) + \<z,\1\> \Big) \cdot y(i) + \Big(\S(y,z) + \S(y,y) + \<y,\1\> \Big) \cdot z(i) +
$$
$$
\Big(\S(y,y) \cdot \<z,\1\> + \S(z,z) \cdot \<y,\1\> + \<y,z\> \cdot \<y+z,\1\>\Big)
$$
By assumption, $\S(y,z) = \<y,\1\> \cdot \<z,\1\> + \<y,z\> = 0$. Note that this also implies $\<y,z\> \cdot \<y+z,\1\> = 0$, implying
$$
\ell(i) =
\Big(\S(z,z) + \<z,\1\> \Big) \cdot y(i) + \Big((S(y,y) + \<y,\1\> \Big) \cdot z(i) +
\Big(\S(y,y) \cdot \<z,\1\> + \S(z,z) \cdot \<y,\1\>\Big)
$$
Consequently, the linear part of $Q$ may be written as
$$
\sum_i \ell(i) x(i) =
$$
$$
\Big(\S(z,z) + \<z,\1\> \Big) \cdot \<x,y\> + \Big((S(y,y) + \<y,\1\> \Big) \cdot \<x,z\> + \Big(\S(y,y) \cdot \<z,\1\> + \S(z,z) \cdot \<y,\1\>\Big) \cdot \<x,\1\>
$$
This means that the non-zero Fourier coefficients of the polynomial $Q = \sum_{i < j} q_{i,j} x(i) x(j) + \sum_i \ell(i) x(i) + c$ lie in the affine subspace $AF_{y,z} = yz + \mbox{Span}\(y,z,\1\)$.

\subsection{Second derivatives of a fixed polynomial of degree $3$}
Let
$$
g(x) = \sum_{i < j < k} a_{i,j,k} x(i) x(j) x(k)
$$
be a polynomial of degree $3$. For directions $y,z \in \FF^N$, consider the second derivative $g_{y,z} = \sum_i v_{y,z}(i) x(i) + c_{y,z}$. We need to show that the probability of the vector $v_{y,z}$ falling in the affine space $AF_{y,z} = yz + \mbox{Span}\(y,z,\1\)$ is exponentially small.

First, some notation. For $1 \le i \le N$, let $G_i$ be a symmetric $N \times N$ matrix over $\FF$ with
$\(G_i\)_{j,k} = \(G_i\)_{k,j} = a_{i,j,k}$ for all $j \not = k$. (Here we think about $\{i,j,k\}$ as an unordered subset of $[N]$.) The diagonal entries of $G_i$ are set to $0$. For future use note the important property $\(G_i\)_{j,k} = \(G_j\)_{i,k} = \(G_k\)_{i,j}$.

These matrices are relevant because they describe the vector $v_{y,z}$.
\lem
\begin{itemize}
\item
$$
v_{y,z}(i) = \mbox{coef}_{x(i)}\(g_{y,z}(x)\) = \<y,G_i z\>
$$
\item
An alternative representation of $v_{y,z}$ will be more convenient for us.
For $z \in \FF^N$, let $G(z) = \sum_{i=1}^N z(i) G_i$. Then
$$
v_{y,z} = G(z) \cdot y
$$
\end{itemize}
\elem
\prf
For the first claim of the lemma, by linearity of the derivative, it suffices to consider the monomial $g(x) = x(i) x(j) x(k)$. This case can be easily verified directly.

For the second claim, note that
$$
(G(z) \cdot y)(l) = \sum_{k=1}^N \(G(z)\)_{k,l} y(k) = \sum_{k=1}^N y(k) \cdot \sum_{i=1}^N z(i) \(G_i\)_{k,l} =  \sum_{k=1}^N y(k) \cdot \sum_{i=1}^N \(G_l\)_{k,i} z(i)= \<y,G_l z\>
$$
\eprf

Consider the event $\{v_{y,z} \in AF_{y,z}\}$. This means $v_{y,z} = yz + u_{y,z}$, for some vector $u_{y,z} \in \mbox{Span}(y,z,\1)$. There are only $8$ possible choices for $u_{y,z}$. For convenience, let us assume, without loss of generality (as can be easily seen from the proof), that $u_{y,z} = y + z + \1$ is the most popular one. By the lemma, the event $\{v_{y,z} = yz + u_{y,z}\}$ is the same as $\{G(z) \cdot y = yz + u_{y,z}\}$. To simplify things some more, let $A_i = G_i + e_i \otimes e_i$, $i=1...N$. That is, $A_i = G_i$ but for $\(A_i\)_{i,i} = 1$. Let $A(z) = \sum_{i=1}^N z(i) A_i$. Note that $A(z)\cdot y = G(z) \cdot y + yz$. Hence $\{G(z) \cdot y = yz + u_{y,z}\}$ is the same as $\{A(z) \cdot y =  u_{y,z} = y + z + \1\}$

We conclude the proof by a technical claim.
\pro
\label{pro:yz-hard-to-get}
Let $\{A_i\}$, $i=1...N$ be a family of symmetric $N\times N$ matrices over $\FF$ with $A_i(k,k) = \delta_{ik}$. Then, for  $y,z$ uniformly at random and independently from $\FN$,
$$
Pr_{y,z} \Big\{\(A(z)\) \cdot y = y + z + \1\Big\} \le \(\frac34\)^N
$$
\epro

The proof of the proposition is based on the claim that the rank of a matrix $A(z)$ is typically large.
\lem
\label{lem:rank-large}
Let matrices $\{A_i\}$ be as in the proposition. Let $C$ be any fixed symmetric $N \times N$ matrix. Then
$$
Pr_z\Big\{\mbox{rank}(A(z) + C) \le k-1\Big\} \le \frac{1}{2^N} \cdot \sum_{i=0}^{k-1} {N \choose i}.
$$
\elem
\prf
Consider a family of ${N \choose k}$ polynomials $f_I$ on $\FN$. These polynomials are indexed by $k$-subsets of $[N]$. For a $k$-subset $I$, let $f_I(z)$ be the determinant of the $I\times I$ minor of $A(z) + C$. Clearly, rank of $A(z) + C$ is smaller than $k$ if and only if $z$ is a joint zero of $\{f_I\}$.

We now claim that the coefficient of $\prod_{i \in I} z_i$ in $f_I(z)$ is $1$. If this is true, $\mbox{deg}(f_I - \prod_{i \in I} z_i) \le k - 1$, and the claim of the lemma will follow from Lemma~\ref{lem:rank-large}.

Let $B(z) = A(z) + C$. Since we are working in characteristic two, the symmetry of $B(z)$ implies that
$$
\det B(z)= \sum_{\sigma \in S_N:~\sigma = \sigma^{-1}} ~~\prod_{i=1}^N B_{i \sigma(i)}(z)=
$$
$$
\sum_{\sigma \in S_N:~\sigma = \sigma^{-1}}
~~\prod_{\{i:\sigma(i)=i\}} \(z_i + C_{i,i}\) \cdot \prod_{\{i: i<\sigma(i)\}}  B_{i \sigma(i)}(z) = \prod_{i \in I}^n z_i + \mbox{lower order terms}.
$$
In the second equality we use the identity $B^2_{i \sigma(i)}(z) = B_{i \sigma(i)}(z)$ in $\FF$.
\eprf

Let $I$ denote the identity $N\times N$ matrix.

Let $p(z) = Pr_y \Big\{A(z) \cdot y = y + z + \1\Big\}$. Clearly $p(z) \le 2^{-\mbox{rank}(A(z)+I)}$. By Lemma~\ref{lem:rank-large},
$$
Pr_{y,z} \Big\{\(A(z)\) \cdot y = y + z + \1\Big\} = \E_z p_z \le \E_z 2^{-\mbox{rank}(A(z)+I)} \le \frac{1}{2^N} \sum_{k=0}^N {N \choose k} 2^{-k} = \(\frac34\)^N
$$
This concludes the proof of the proposition, and of Theorem~\ref{thm:non-approx}.

\section{Acknowledgements}
We are grateful to Ben Green and Terence Tao
for sending us their paper. We would also like to thank them and
Emanuele Viola for bringing the argument of Alon and Beigel to our
attention


\begin{thebibliography}{99}

\bibitem{ab}
N. Alon, R. Beigel {\sl Lower Bounds for Approximations by Low Degree Polynomials over $Z_m$},
{SCT}: Annual Conference on Structure in Complexity Theory, 2001.

\bibitem{bv}
A. Bogdanov, E. Viola {\sl Pseudorandom bits for polynomials via the Gowers norm}, FOCS'07.

\bibitem{gowers}
W. T. Gowers, {\sl A new proof of Szemeredi's theorem}, GAFA
Vol. 11(2001), pp. 465-588.

\bibitem{gt}
B. Green, T. Tao, {\sl An inverse theorem for the Gowers $U^3$ norm},
Proc. Edinburgh Math. Soc., to appear.

\bibitem{gt07}
B. Green, T.Tao, {\sl The distribution of polynomials over finite fields,
with applications to the Gowers norms}, preprint, 2007.

\bibitem{MS}
J. MacWilliams and N. J. A. Sloane, {\bf The Theory of Error
Correcting Codes}, Amsterdam, North-Holland, 1977.

\bibitem{sam}
A. Samorodnitsky, {\sl Low degree tests at large distances}, STOC '07.

\bibitem{tao-lec}
T. Tao,  {\sl Structure and randomness in combinatorics}, FOCS '07.

\bibitem{vw}
E. Viola, A. Wigderson, {\sl Norms, XOR lemmas, and lower bounds for $GF(2)$
polynomials and multiparty protocols}, CCC '07.

\end{thebibliography}
\end{document}